\documentclass[11pt]{article}
\makeatletter
\def\section{\@startsection{section}{1}{\z@}{-3.5ex plus -1ex minus -2.ex}
{2.3ex plus .2ex}{\Large\bf}}

\def\subsection{\@startsection{subsection}{2}{\z@}{-3.25ex plus
 -1ex minus -2.ex}
{1.5ex plus .2ex}{\bf}}

\def\vsn{\vskip 1pc \noindent}
\def\f{\newline}
\def\e{\varepsilon}
\def\a{\alpha}

\def\comp{{\rm comp}}

\def\rr{{\bf R}}
\def\nn{{\bf N}}

\def\vN{\vec{N}}
\def\vM{\vec{M}}

\newcommand{\be} {\begin{equation}}
\newcommand{\ee} {\end{equation}}
\newcommand{\bd} {\begin{displaymath}}
\newcommand{\ed} {\end{displaymath}}

\newcommand{\bq}{\begin{eqnarray}}
\newcommand{\eq}{\end{eqnarray}}
\newcommand{\bqn}{\begin{eqnarray*}}
\newcommand{\eqn}{\end{eqnarray*}}
\newcommand{\ba}[1]{\begin{array}{#1}}
\newcommand{\eqa}{\end{array}}
\def\qed{
   \\[-4ex]
  \hbox to \hsize{\hfill \vrule height 1.6ex width 1.5ex
  depth -.1ex}}

\begin{document}

\bibliographystyle{alpha}

\begin{center} {\Large {\bf 
Efficient  finite-dimensional solution of initial value problems in infinite-dimensional Banach spaces 
 } } 
\end{center}
\footnotetext[1]{
The author was partially  supported by the Polish
NCN grant - decision No. DEC-2017/25/B/ST1/00945    and by the Polish   Ministry of Science  and Higher Education
}              
  
\vspace{1mm}
\begin{center}
{\large Boles\l aw Kacewicz \footnotemark[1] \footnotemark[3] $\;\;\;$ Pawe\l $~$ Przyby\l owicz}\footnotemark[2] \footnotemark[3] \\
\vspace{1mm}
{\it 
AGH University of Science and Technology\\
 Faculty of Applied Mathematics\\
Al. Mickiewicza 30, 
30-059 Krakow, Poland \\
}

\end{center}

\footnotetext[2]{ 
The author was partially supported by the Faculty of Applied Mathematics AGH UST dean grant No.15.11.420.038/1 for PhD students and young researchers within subsidy of Ministry of Science and Higher Education
}

\footnotetext[3]{ 
E-mail: 
B. Kacewicz {\tt kacewicz@agh.edu.pl} (corresponding author)\\
$~$ \hspace{16mm} P. Przyby\l owicz {\tt pprzybyl@agh.edu.pl} 
}

\thispagestyle{empty}

\begin{center}
21.02.2018  $\;\;\;$     
\end{center}

\vsn
\begin{center} {\bf{\Large Abstract}} \end{center}
\noindent
We deal with the approximate solution of initial value problems in infinite-dimensional Banach spaces with a Schauder basis. We only allow
 finite-dimensional algorithms acting in the spaces $\rr^N$, with varying $N$.  The error of such  algorithms depends on 
two parameters: the truncation parameters $N$ and a discretization parameter $n$. 
For a class of $C^r$ right-hand side functions, we define an algorithm with varying $N$, based on possibly non-uniform mesh, and
we analyse its error and cost. For constant $N$, we
show a matching (up to a constant)  lower 
bound on the error of
any algorithm in terms  of  $N$ and $n$, as $N,n\to \infty$. We stress that 
in the standard error analysis  the dimension $N$ is fixed, and  the dependence on $N$ 
is usually hidden in  error coefficient.  For a certain model of cost,  for many cases of interest, we show tight 
(up to a constant)  upper and lower bounds on the minimal cost
of computing an $\e$-approximation to the solution (the $\e$-complexity of the problem). The results are
illustrated by an example of the initial value problem in the weighted $\ell_p$ space ($1\leq p<\infty$).
\f

\vspace{1cm}

\newpage
\noindent
{\Large \section{ Introduction }}
\noindent
Let  $(E, \|\cdot\|)$ be an infinite-dimensional  Banach space over $\rr$ with a Schauder basis. 
We study the solution of an initial value problem
\be
z'(t)=f(z(t)),\;\;\; t\in [a,b],  \;\;\;\; z(a)=\eta,
\label{rnie}
\ee
where  $a<b$, $\eta\in E$ and $f: E\to E$ is a Lipschitz function in $E$.
The Lipschitz condition implies the existence and uniqueness of a solution   $z:[a,b]\to E$, see e.g. \cite{Cartan}, \cite{deiml} or \cite{LusSob}.
We aim at approximating the solution $z$ in $[a,b]$. 
\f
Infinite countable systems of the form (\ref{rnie}) have been investigated  for many years. They often appear in various applications inspired 
by physical, chemical or mechanical problems, see, for example, \cite{Bel} ,  \cite{BelAdo}, \cite{deiml}, \cite{Lewis}  or  \cite{Temam}. 
Many basic results have already been surveyed  in \cite{deiml}. 
According to \cite{deiml}, one  can distinguish between two main approaches to countable systems: a direct approach, where 
we look for a sequence $z(t)$ satisfying the sequence of equations (\ref{rnie}), and a Banach space approach, where 
$f$ acts in a Banach space, and the solution $z$ is a Banach space valued function. 
\f
In this paper we consider computational aspects of (\ref{rnie}).
In contrast to   theoretical properties of infinite systems,  much less is known about efficient approximation of the solutions, 
see e.g. \cite{Brzych}, \cite{deiml}. 
The authors  of most papers  concentrate  on basic Galerkin-type devices  that allow us to truncate the infinite system to a finite-dimensional one.
An extensive complexity    analysis of problem (\ref{rnie}) in  a  Banach space 
in the deterministic and randomized settings   is recently presented in \cite{Hein1}, \cite{Hein2} and \cite{Hein3}. 
The authors  assume  that computations in the underlying Banach space are allowed.  This means, in particular,  that the evaluations 
of the Banach space valued right-hand side function $f$ and its partial derivatives can be performed with  the unit cost. Complexity upper bounds are obtained by 
a complex multilevel projection algorithm. 
\f
It is well known that in the finite-dimensional case, for $\rr^N$-valued functions $f$ with finite and fixed $N$,  there is a vast literature devoted to   optimal approximation of the solution of 
  (\ref{rnie}), see, for example,  \cite{daun} , \cite{bk1},
\cite{bk2},  \cite{bkpp1},  or many other papers.   
\f  
Motivated by computer applications,  we 
 restrict  ourselves in this paper to algorithms for (\ref{rnie}) that are only based on finite-dimensional computations. No operations performed in $E$ are  allowed.
In particular, we do not admit computation of $f(y)$ for $y\in E$.  
Thus, we consider a different computational model compared to 
that in \cite{Hein1}.  
We wish to study the quality of such  finite-dimensional solution of (\ref{rnie}).
\f
Our approach is different from that in \cite{Hein1}, \cite{Hein2} and \cite{Hein3}.  We assume that the space $E$ has a Schauder basis. Since most important spaces appearing in applications, such as
$\ell_p$ or $L_p$ for $1\leq p<\infty$, have  Schauder bases, this is not a restrictive assumption in practice. 
The class of problems under consideration and the class of algorithms are defined in terms of that basis. 
Our main results  are as follows:
\vsn
$\bullet$ For a class of $C^r$ functions $f: E\to E$ we define an algorithm $\phi^*_{n,\vN}$ based on possibly non-uniform mesh
(with $n+1$ points) 
and restricted, in each time step, to finite-dimensional computations with  varying dimensions, represented by the vector $\vN$.
We show an upper bound on the error  of $\phi^*_{n,\vN}$ expressed in the terms of the truncation vector $\vN$  and 
the step sizes. 
In contrast to the usual analysis in the finite-dimensional case, the parameter $N$ is now  not a constant number, which may be hidden in 
error coefficient,  but it tends to infinity. 
This requires somewhat different analysis
including the tractability questions, see \cite{NW}.
\f
$\bullet$ For constant dimensions equal to $N$, we bound from below the error of any algorithm $\phi_{n,\vN}$ for solving (\ref{rnie}). 
The bound shows that the algorithm $\phi^*_{n,\vN}$ is error optimal (up to a constant)  as $n,N\to \infty$.
\f
$\bullet$ Based on two-sided error bounds, for $\e>0$ we discuss upper and lower bounds on the minimal cost of computing an $\e$-approximation to the solution of (\ref{rnie})
(the $\e$-complexity of the problem).  To consult  a general notion of the $\e$-complexity,  see \cite{TWW}.
\f
$\bullet$ We illustrate the results by an example of a countable system of equations in $\ell_p$, $1\leq p< \infty$, embedding it to the Banach space setting 
with a weighted $\ell_p$ space.
\vsn
The paper is organized as follows.
In Section 2 we present basic notions and definitions, and we define the model of computation. In Section 3 we define the algorithm $\phi^*_{n,\vN}$ and prove an upper error bound in
Theorem 1. Theorem 2 shows a lower error bound for any algorithm $\phi_{n,\vN}$ based on constant truncation parameters.
In Propositions 1 nad 2 we discuss the resulting $\e$-complexity bounds for the problem (\ref{rnie}).
 Section 4 contains an example of a countable system, to which we apply the results described in Section 3. We show how to select $N$ and $n$ to get the error at most $\e$,
and we establish  the cost of computing the $\e$-approximation.
 In Section 5 we recall, for convenience of the reader, basic facts used in the paper  about integration, differentiation and interpolation  in a Banach space with Schauder basis. 
\vsn
\noindent
{\Large \section{ Preliminaries }}
\noindent
\vsn
 Let $\{e_1,e_2,\ldots \; \}$ with  $\|e_j\|=1$ be a Schauder basis in $E$.  Let $f(y)=\sum\limits_{j= 1}^{\infty}  f^j(y)e_j$ for $y\in E$. 
For  $k\in \nn$,  let $P_k: E\to E$ be the projection operator, i.e.,  for $z\in E$,  $z=\sum\limits_{j= 1}^{\infty} z^je_j$ we have
$P_kz:= \sum\limits_{j=1}^k z^je_j$.
  The operator $P_k$ is linear and bounded, 
with  $\sup\limits_{k} \|P_k\|=P<\infty$, see \cite{LinTza} p. 1--2. The number $P$ is called the basis constant of $\{e_1,e_2,\ldots \;\}$. Since $|z^k|= \|(P_k-P_{k-1})z\|$ ($P_{0}=0$),  it holds 
$|z^k|\leq 2P\|z\|$, for $k=1,2,\ldots \;\;$.
\vsn
{\it The class of problems}
\vsn
Let $r$ be a nonnegative integer. Let $L$, $M$, $D$ be positive numbers, and $\Gamma=\{\gamma(k)\}_{k=1}^{\infty}$ and   $\Delta=\{\delta(k)\}_{k=1}^{\infty}$    positive, nonincreasing,  convergent to zero sequences. 
We shall consider a class $F_r=F_r(L,M,D, P,\Gamma, \Delta)$ of pairs $(f,\eta)$ defined  by the following conditions (A1)--(A5).   
\vsn
(A1)$\;\;$ $ \| \eta-P_k\eta\|\leq \gamma(k)$ for $k\in \nn$,
\f
(A2)$\;\;$  $\|f(y)-f(\bar y)\|\leq L\| y-\bar y\|$, for $y,\bar y\in E$,
 \f
(A3)$\;\;$  $\|f(\eta)\|\leq M$.
\vsn
 Let $R=R(L,M,P,a,b)$  be a number, existence of which is shown in Lemma 1 below,  and let $K=K(\eta,R)=\{y\in E: \|y-\eta\| \leq R\}$.
In addition to (A1)--(A3), we assume  that
\vsn
(A4)$\;\;$ $f\in C^r(K)$ (where the derivatives are meant in the  Fr\'echet sense) and
 $$\sup\limits_{y\in K} \|f^{(l)}(y)\| \leq D, \;\;\;\; l=1,2,\ldots,r. $$ 
(In the last inequality,  $\|\cdot\|$ means the norm of a bounded $l$-linear operator in $E^l$, defined by the norm in $E$.) 
We note  that for $l=0$ we have from (A2) and (A3) the bound  
$$
\sup\limits_{y\in K} \|f(y)\| \leq M+LR.
$$
(A5)$\;\;$  $\sup\limits_{y\in K} \left\|f(y)-P_kf(y)\right\| \leq \delta(k),\;\;\; k\in \nn.$  
\vsn
The parameter $P$ of the space $E$, as well as the parameters $L$, $M$, $D$, $\{\gamma(k)\}_{k=1}^\infty$, $\{\delta(k)\}_{k=1}^\infty$ 
of the class $F_r$ are unknown, and they cannot be used by an algorithm. The numbers $a$, $b$, $r$ are known.
\vsn
{\it The class of algorithms}
\vsn
To approximate $z$,  we shall only allow  Galerkin-type algorithms  that base on finite-dimensional computations. 
 Let $n\in \nn$ be a discretization parameter,  and let $\{\alpha(n)\}_{n=1}^\infty$ be a nonincreasing sequence convergent to $0$
as $n\to \infty$.  We shall consider  a family of partitions of $[a,b]$ given by points
$a=t_0^n<t_1^n<\ldots< t_n^n=b$ such that  
\be
 \max_{ k=0,1,\ldots, n-1  } (t_{k+1}^n-t_k^n) \leq \alpha(n), \; n=1,2,\ldots \, .  
\label{part}
\ee
(Obviously, it must hold $\alpha(n)\geq (b-a)/n$.) 
In what follows, we shall omit in the notation the superscript $n$.  Furthermore,   to keep the notation legible, 
we will not indicate the dependence of information and an algorithm on $\{t_k\}_{k=0}^n$. 
\f
Let $N_{-1}\in \nn$ and $N_k, M_k\in \nn$,  $k=0,1,\ldots, n-1$. For a given partition $\{t_k\}_{k=0}^{n}$ and given numbers 
$\{N_k\}_{k=-1}^{n-1}$ and $\{M_k\}_{k=0}^{n-1}$, 
we shall allow algorithms  based on the (approximate) successive solution
of finite-dimensional local problems in $[t_k,t_{k+1}]$, $k=0,1,\ldots,n-1$.  .
Let  $\bar f_{k}=P_{N_k} f$, and let $\bar z_{k}:[t_k, t_{k+1}]\to E$ be the solution of the local problem 
\be
\bar z_{k}'(t)=\bar f_{k}(\bar z_{k}(t)),\;\; \bar z_{k}(t_k)=P_{M_k} y_k,\;\;\; t\in [t_k,t_{k+1}],
\label{rnielocal1}
\ee
where $y_k$ is a given point in $E$. This is a finite-dimensional problem defined by truncation parameters $N_k$ 
(which describes  the number of components of $f$ that are considered) and $M_k$ (which describes the number 
of components of the arguments  taken into account).
Note that  $\bar f_{k}$ is a Lipschitz function in $E$   with the (uniform) constant $PL$. 
An algorithm  successively computes $y_k$ and  approximations $l_k$ to  $\bar z_k$ in $[t_k,t_{k+1}]$.
  The approximation to the solution 
$z$ of (\ref{rnie})  in $[a,b]$ is a spline function $l:[a,b]\to E$ composed of $l_k$. 
\f
Consider information about $f$  that is allowed in the computation of $l$. 
The function $f$ can  only be accessed through  the components $\bar f_{k}^j$ of $\bar f_{k}$,  $j=1,2,\ldots,N_k$.
Available   information is given by evaluations  $\bar f_{k}^j(P_{M_k} y)$,  or evaluations  of
 partial derivatives of  $\bar f_{k}^j(P_{M_k} y)$ at $P_{M_k} y$  (up to order $r$),  for some component $j$,  at some  
information  points $y$.  We assume that
the number of  information points $s$  is proportional to the number of subintervals, that is,  there is $\hat K$ such that
\be
s\leq \hat Kn,
\label{hatK}
\ee
$n=1,2,\ldots \;$.
For example,  the standard explicit Runge-Kutta methods in $\rr^N$ of order $p$  are based, at each time interval
$[t_k, t_{k+1}]$, on   the constant number of $p$ function evaluations.   For what concerns the initial condition $\eta$, we assume 
to have access to $P_{N_{-1}}\eta$ for any $N_{-1}\in \nn$. We assume that information is adaptive in the following sense.
We allow successive adaptive selection of  the information points,  indices  $j$ 
of the components of $\bar f_{k}$,  and orders  of partial derivatives to be evaluated.
This means that  these elements can be computed based on information computed so far. 
In this paper,  the sequences  $\{t_k\}_{k=0}^{n}$,  $\{N_k\}_{k=-1}^{n-1}$ and $\{M_k\}_{k=0}^{n-1}$ are given in advance.
\f
Infinite-dimensional 'computation' is not allowed;  for example,   computing $f(y)\in E$ for  $y\in E$ is in general not possible.
\f
Let $\vN = [N_{-1},N_0,\ldots, N_{n-1}]$ and $\vM = [M_0,M_1,\ldots, M_{n-1}]$.
Information computed as described above in the interval $[a,b]$ for $f$ and $\eta$  will be denoted by ${\cal N} _{n,\vN, \vM}(f,\eta)$. 
By an  algorithm $\phi_{n,\vN, \vM}$ we mean  a mapping that assignes to the vector ${\cal N} _{n,\vN, \vM}(f,\eta)$
the function $l$ described above, $l=\phi_{n,\vN, \vM}\left( {\cal N} _{n,\vN, \vM}(f,\eta)\right)$.
The (worst case)  error of an algorithm  $\phi_{n,\vN, \vM}$ with information  ${\cal N} _{n,\vN, \vM}$  in the class $F_r$ is defined by
\be
e(\phi_{n,\vN, \vM}, {\cal N} _{n,\vN, \vM} ,F_r)=\sup\limits_{(f, \eta)\in F_r} \, \sup\limits_{t\in [a,b]}\, \| z(t)-l(t)\|.
\label{error}
\ee
Let us consider the cost of computing information ${\cal N} _{n,\vN, \vM}(f,\eta)$. 
For each $j$, we assume that the cost of computing 
the value of the function $\bar f_{k}^j$ or its partial derivative at  $P_{M_k} y$ is  $c(M_k)$, where $c$ is a nondecreasing function. 
That is, the cost of computing  these real-valued functions is determined by  the number of variables $M_k$.
\f
The  number of such scalar evaluations at each time step $[t_k,t_{k+1}]$ 
depends on  particular information; we denote this number by  $\ell(N_k,M_k)$.  For instance, if we only compute at each time step a single value
$\bar f_{k}(P_{M_k} y)$ for some $y$, then $\ell(N_k,M_k)=N_k$. If we compute at  $P_{M_k} y$
all partial derivatives up to order $r$ of each component of $\bar f_{k}$, then
$\ell(N_k,M_k)=\Theta(N_k M_k^r)$. 
\f
We assume to have access to $P_{N{-1}}\eta$ for any $N_{-1}$ with no cost.
\f
The total cost of computing information is thus 
$$
\sum\limits_{k=0}^{n-1}  c(M_k)\ell(N_k,M_k).
$$
For a given $\e>0$, by the $\e$-complexity $\comp(\e)$ of the problem (\ref{rnie}), we mean the minimal cost of computing an $\e$-approximation.
More precisely,
$$
\comp(\e)=\inf \left\{  \sum\limits_{k=0}^{n-1}  c(M_k)\ell(N_k,M_k)   : \; n, \vN, \vM \mbox{ are such that } \exists \, {\cal N} _{n,\vN, \vM}, \,
 \phi_{n,\vN, \vM} \right.
$$
\be
 \mbox{ with }  \left. e(\phi_{n,\vN,\vM}, {\cal N} _{n,\vN, \vM},  F_r)\leq \e \right\}      .
\label{comp}
\ee
The $\e$-complexity measures an intrinsic difficulty of solving the problem (\ref{rnie}) by finite-dimensional computation.
We shall establish in this paper  bounds on $\comp(\e)$.
\f
Unless otherwise stated, all coefficients that appear in this paper will only depend on $L$, $M$, $P$, $D$, $r$, $a$ and $b$.
\vsn
\noindent
{\Large \section{ Upper error and complexity bounds }}
\noindent
\subsection{The variable dimension algorithm}
\noindent 
Let $h_k=t_{k+1}-t_k$.  
We associate with each subinterval $[t_k,t_{k+1}]$ a number (dimension) $N_k$, $k=0,1,\ldots,n-1$.  For $\vN=[N_{-1}, N_0, \ldots, N_{n-1}]$, 
we  define an  algorithm $\phi^*_{n,\vN}$ for solving (\ref{rnie}).  Let $\bar y_0=P_{N_{-1}}\eta$.
 Let $\bar f_k=P_{N_k} f$ and, for a given $\bar y_k\in E$, consider the local problem
\be
\bar z_{k}'(t)=\bar f_{k}(\bar z_{k}(t)),\;\; \bar z_{k}(t_k)=\bar y_k,\;\;\; t\in [t_k,t_{k+1}]. 
\label{rnielocal1a}
\ee
The following general idea of approximating the solution $z$ of (\ref{rnie}) has been used several times in various contexts,  see e.g. \cite{daun}, \cite{Hein1},  \cite{bkpp1}.  
 Let $r\geq 1$.     We define a function $\bar  l_{k, r}$ in $[t_k,t_{k+1}]$ as follows. 
Let $\bar l_{k,0}(t) \equiv \bar y_{k}$. For $s\geq 0$ and a given function $\bar l_{k,s}$, we define the Lagrange interpolation 
polynomial of degree at most $s$  by
\be
\bar q_{k,s}(t)= \sum\limits_{p=0}^s\prod\limits_{l=0,\,  l\ne p}^s \, \frac{t -\xi_{k,l}}{\xi_{k,p}-\xi_{k,l} }  \bar f_k(\bar l_{k,s}(\xi_{k,p})) ,
\label{polyn}
\ee
where $\xi_{k,p}= t_k+ ph_k/s$,  $p=0,1,\ldots, s$  (with $ \prod\limits_{l=0,l\ne p}^0=1$). We  define a polynomial 
\be
\bar l_{k,s+1}(t)=\bar y_{k} + \int\limits_{t_k}^t  \bar q_{k,s}(\xi)\, d\xi, \;\;\;\; t\in [t_k,t_{k+1}].
\label{hatl}
\ee
We repeat (\ref{polyn}) and  (\ref{hatl}) for $s=0,1,\ldots, r-1$ to get  a final polynomial $\bar l_{k,r}$, and we set
 $\bar y_{k+1}= \bar l_{k,r} (t_{k+1})$.  After passing through all the intervals $[t_k,t_{k+1}]$, $k=0,1,\ldots, n-1$, 
we get an approximation to $z$ in $[a,b]$ defined as a piecewise polynomial continuous function $\bar l_n(t):=\bar l_{k,r}(t)$,
if $t\in [t_k,t_{k+1}]$.  
\f
   For $r=0$ we define $\bar l_n$ as we did above for $r=1$, i.e.,
$$\bar l_n(t)= \bar y_{k} +(t-t_k)\bar f_k(\bar y_{k}), \;  \mbox{ if }\; t\in [t_k,t_{k+1}],$$ 
and  we set as above $\bar y_{k+1}= \bar l_{n} (t_{k+1})$.
\f
Note that the computation of $\bar l_n$ only involves  finite-dimensional operations. We see that 
the computations are determined by the vector $\vN=[N_{-1},N_0,\ldots, N_{n-1}]$, since  the dimensions $M_k$  are 
given by $M_k=\max\limits_{j=-1,0,\ldots,k} N_j$ for $k=0,1,\ldots, n-1$.
\f
We denote the information about $f$ used above to construct $\bar l_n$ by 
${\cal N}^*_{n,\vN}(f,\eta)$.  It is easy to see that  it consists of $O(r^2n)$ evaluations of finite-dimensional truncations 
of $f$ at finite-dimensional truncations of some points in $E$. 
We define an algorithm $\phi^*_{n,\vN}$ that approximates $z$ by  
\be 
\phi^*_{n,\vN} \left( {\cal N}^*_{n,\vN}(f,\eta)\right)(t) = \bar l_n (t), \;\; t\in [a,b].
\label{algorithm}
\ee
\noindent
\subsection{Upper error bound}
\noindent 
We show in this section an upper bound on the error of $\phi^*_{n,\vN}$.  We start with a lemma that 
assures that the solutions of  (\ref{rnie}), (\ref{rnielocal1a})  stay in a certain ball $K(\eta,R)$ with radius $R$ that only depends on
the parameters appearing in assumptions (A1), (A2) and (A3).
\vsn
{\bf Lemma 1}$\;\;$ {\it There exist  $R=R(L,M,\gamma(1), P,a,b)$ such that:
\f
{\it (a)} $\;\;$  $z(t)\in K(\eta,R) \mbox{ for } t\in [a,b]$
\f
and 
\f
{\it (b)} $\;\;$ for any $\{\alpha(n)\}_{n=1}^\infty$  there is $\hat n$
such that for any $n\geq \hat n$,  
any  $\{t_k\}_{k=0}^{n}$ satisfying (\ref{part}), any $\vN$ 
and any $f$ satisfying (A1), (A2) and (A3) it holds 
$$
\bar z_{k}(t)\in  K(\eta, R) \mbox{ for } t\in [t_k,t_{k+1}], \;\; k=0,1,\ldots, n-1 .  
$$
}
{\bf Proof}$\;\;$ Note that the boundedness of $\|z(t)- \eta\|\leq \tilde C_1$ 
by some constant $\tilde C_1=\tilde C_1(L,M, P,a,b)$ immediately follows from the Gronwall inequality,  (A2) and (A3).
  We  show a bound on $\| \bar z_{k}(t) - \eta\|$.  
For $r=0$  the algorithm is defined by the same expression as for $r=1$, so that we can consider 
formulas  (\ref{polyn}) and (\ref{hatl})  with $r\geq 1$. We have  
$$
\bar l_{k,s+1}(t) - \bar y_{k}= 
$$
$$\int\limits_{t_k}^t \left(\sum\limits_{p=0}^s
\prod\limits_{l=0,\,  l\ne p}^s \, \frac{\xi -\xi_{k,l}}{\xi_{k,p}-\xi_{k,l} }  \left(\bar f_k(\bar l_{k,s}(\xi_{k,p})) - \bar f_k(\bar y_{k}) \right) 
 +   \bar f_k(\bar y_{k})    \right)\, d\xi. 
$$
Since, $\bar f_k$ is the Lipschitz function  with the constant $PL$, we have that
\be
\| \bar l_{k,s+1}(t) - \bar y_{k}\| \leq h_k PL C_r \sup\limits_{t\in [t_k,t_{k+1}]}  \| \bar l_{k,s}(t) - \bar y_{k} \| +h_k \| \bar f_k(\bar y_{k}) \|, 
\label{pom1}
\ee
$t\in [t_k,t_{k+1}]$, $s=0,1,\ldots, r-1$,  where $C_r$ is an upper bound (dependent only on $r$) on 
$$
\sup\limits_{\xi\in [t_k,t_{k+1}]}\, \sum\limits_{p=0}^s \prod\limits_{l=0,\,  l\ne p}^s \, \left| \frac{\xi -\xi_{k,l}}{\xi_{k,p}-\xi_{k,l} } \right|.
$$
Taking the $\sup\limits_{t\in [t_k,t_{k+1}]}$ in the left hand side of (\ref{pom1}), and solving the resulting recurrence inequality, we get for 
$n$ such that $2\alpha(n) PLC_r\leq 1$ the bound
\be
 \sup\limits_{t\in [t_k,t_{k+1}]}  \| \bar l_{k,s}(t) - \bar y_{k} \| \leq 2h_k \| \bar f_k(\bar y_{k}) \|\leq 2h_k P\| f(\bar y_k)\| ,
\label{pom2}
\ee
$s=0,1,\ldots,r$.    Setting   $s=r$,  $t=t_{k+1}$,  in the left-hand side of  (\ref{pom2}) we also have
\be
\| \bar y_{k+1} - \bar y_{k} \| \leq  2h_kP\| f(\bar y_k)\| .  
\label{pom3}
\ee
 We now  bound $E_k= \|\bar y_{k} - \bar y_{0}\|$ ($E_0=0$). Since 
$ \|f (\bar y_{k}) \| \leq \| f(\bar y_{0}) \| + L E_k$, we get from (\ref{pom3}) that
$$
E_{k+1}\leq (1+2h_kPL) E_k + 2h_k  P\| f(\bar y_{0}) \|, \;\;\; k=0,1,\ldots, n-1.
$$
By solving the recurrence inequality, remembering that $\sum\limits_{j=0}^k h_j = t_{k+1}-t_0\leq b-a$,  we get 
\be
\|\bar y_{k} - \bar y_{0}\| \leq C\, \|  f(\bar y_{0}) \|, \;\; k=0,1,\ldots, n-1, 
\label{pom4}
\ee
for some constant $C$ only dependent on $L,P,a,b$ and sufficiently large $n$.
\f
We now estimate $\|\bar z_{k}(t)-\eta\|$, $t\in [t_k,t_{k+1}]$. We have from (\ref{rnielocal1a}) that
$$
\bar z_{k}(t) = \bar y_k+ \int\limits_{t_k}^t \bar f_k(\bar z_{k}(\xi))\, d\xi.
$$
Subtracting from both sides $\bar y_0$ and  applying the Lipschitz condition for $\bar f_k$, we get
$$
\| \bar z_{k}(t) -\bar y_0\| \leq \| \bar y_{k} - \bar y_0 \|   + PL  \int\limits_{t_k}^t \| \bar z_{k}(\xi) -\bar y_0\| \, d\xi +h_k P\|f(\bar y_0) \|, \; t\in [t_k,t_{k+1}].
$$
By  Gronwall's lemma,  we get
$$
\| \bar z_{k}(t) - \bar y_0 \| \leq \exp(h_kPL) \left( \| \bar y_{k} - \bar y_{0}\| +h_k P \|  f(\bar y_0)\| \right), 
$$
which yields according to  (\ref{pom4})  that
\be
\| \bar z_{k}(t) -\bar y_0 \| \leq  C  \| f(\bar y_0)\| , \;\; t\in [t_k,t_{k+1}],\;\; k=0,1,\ldots, n-1, 
\label{pom5}
\ee
for some constant $C$ (different from that in (\ref{pom4})) and sufficiently large $n$. 
\f
Since  $\bar y_0=P_{N_{-1}}\eta$,  it follows from (\ref{pom5}) and (A1), (A2) and (A3) that
$$
\| \bar z_{k}(t) -\bar y_0 \| \leq  C \left( L\gamma(N_{-1})  + M \right)\leq   C \left( L\gamma(1)  + M \right).
$$
Finally, we get  for $n$ sufficiently large, $t\in [t_k,t_{k+1}]$ and $k=0,1,\ldots, n-1$ that
\be
\| \bar z_{k}(t) -\eta\| \leq \| \bar z_{k}(t) -\bar y_0\| +  \| P_{N_{-1}}\eta  -\eta\| \leq  \tilde C_2 ,
\label{pom6}
\ee
for some constant $\tilde C_2$ which only depends on the parameters appearing in assumptions (A1), (A2) and (A3),  $P$, $a$ and $b$.  
To complete the proof,  in the statement of the lemma  we  take $R=\max\{\tilde C_1,  \tilde C_2\}$.
 \qed
\vsn
The following theorem gives us an upper bound on the error of the algorithm $\phi^*_{n,\vN}$. 
In the proof,  we  have to pay attention to  the independence of constants appearing in the bounds of $N_{-1},N_0,\ldots, N_{n-1}$. 
\vsn
{\bf Theorem 1}$\;\;$ {\it  There exist $C$ such that for any $\{\alpha(n)\}_{n=1}^\infty$ there is $\hat n$
such that for  $n\geq \hat n$, any $\{t_k\}_{k=0}^n$  satisfying (\ref{part}), any  $\vN=[N_{-1},N_0,\ldots, N_{n-1}]$,
 any $(f,\eta)\in F_r$ it holds 
\be
\sup\limits_{t\in [t_k,t_{k+1}]}\| z(t)-\bar l_{n}(t)\| \leq C \left( \gamma(N_{-1}) +\sum\limits_{j=0}^{k} h_j \delta(N_j) +    
\sum\limits_{j=0}^{k} h_j^{\max\{r,1\}+1 } \right) ,
\label{thm1-1}
\ee
$k=0,1,\ldots, n-1$.
\f
Consequently, 
\be
e(\phi^*_{n,\vN}, \mathcal{N}^*_{n,\vN}, F_r) \leq C \left( \gamma(N_{-1}) +\sum\limits_{j=0}^{n-1} h_j \delta(N_j) +    
\sum\limits_{j=0}^{n-1} h_j^{\max\{r,1\}+1 } \right) .
\label{thm1}
\ee
Note that the last term in (\ref{thm1}) can be bounded by $\alpha(n)^{\max\{r,1\} } (b-a)$.
}
\vsn
{\bf Proof}$\;\;$ 
  Let $(f,\eta)\in F_r$.  We need to estimate for $t\in [a,b]$ 
$$
\Bigl\| z(t)-\phi^*_{n,\vN}\Bigl(\mathcal{N}^*_{n,\vN}(f,\eta)\Bigr)(t)\Bigl\|.
$$
 We first show that for $r\geq 1$, $k=0,1,\ldots,n-1$ 
\begin{equation}
\label{Lagr_est_lr}
	\sup\limits_{\xi\in [t_k,t_{k+1}]}\|\bar z_{k}(\xi)-\bar l_{k,r}(\xi)\| \leq C_1h_k^{r+1}, 
\end{equation}
for some constant $C_1$ only dependent on the parameters 
of the class $F_r$, and sufficiently large $n$. 
\f
For the knots  $\xi_{k,p}$   given in (\ref{polyn}),  we  let  $\bar w_{k,s}$ be the Lagrange interpolation polynomial 
for $\bar z'_{k}$ defined in a similar way as $\bar q_{k,s}$,
\begin{equation}
	\bar w_{k,s}(\xi)=\sum\limits_{p=0}^{s}\prod_{l=0,l\neq p}^{s}\frac{\xi-\xi_{k,l}}{\xi_{k,p}-\xi_{k,l}} \bar z'_{k}(\xi_{k,p}).
\end{equation}
Note that $\bar w_{k,0}(\xi)=\bar z'_{k}(t_k)=\bar f_k(\bar y_{k})=\bar q_{k,0}(\xi)$.
From   (\ref{hatl})  we have  for $t \in [t_k,t_{k+1}]$ and $s=0,1,\ldots,r-1$  that
$$
	\bar z_{k}(t)-\bar l_{k,s+1}(t)=\int\limits_{t_k}^t \Bigl(\bar z'_{k}(\xi)- \bar w_{k,s}(\xi) \Bigr)d\xi+
\int\limits_{t_k}^t \Bigl(\bar w_{k,s} (\xi)-\bar q_{k,s}(\xi)\Bigr)d\xi,
$$
 which yields 
\be
	\left\| \bar z_{k}(t)-\bar l_{k,s+1}(t)\right\| \leq \int\limits_{t_k}^t \left\| \bar z'_{k}(\xi)- \bar w_{k,s}(\xi) \right\|\,d\xi+
\int\limits_{t_k}^t \left\| \bar w_{k,s} (\xi)-\bar q_{k,s}(\xi)\right\|\, d\xi .
\label{aux1}
\ee
We now bound both terms in the right-hand side of (\ref{aux1}).
Using the integral form of the Lagrange interpolation remainder formula (\ref{interp-remain}), 
 we get for  $\xi\in [t_k,t_{k+1}]$ and $s=0,1,\ldots,r-1$ 
\begin{equation}
\label{est_lagr}	
	\|\bar z'_{k}(\xi)-    \bar w_{k,s}(\xi) \|\leq h_k^{s+1} /(s+1)!\, \sup\limits_{t\in [t_k,t_{k+1}]} \|\bar z^{(s+2)}_{k}(t)\|.
\end{equation}
The Fr\'echet derivative $\bar z^{(s)}_{k}(t)$ for $s=1,2,\ldots, r+1$ can be expressed in the well known way as a sum of multilinear expressions
involving the Fr\'echet derivatives of $\bar f_k= P_{N_k} f$ of order $0,1,\ldots,r$, evaluated at $\bar z_{k}(t)$. Since, 
by  Lemma 1, $\bar z_{k}(t)$ lies in the ball $K$ for sufficiently large $n$, from the assumption (A4) we get for any $N_k$ that
$$
\|\bar z^{(s)}_{k}(t)\| \leq \hat C_2, \;\;\; \mbox{ for } t\in [t_k,t_{k+1}],\;\; k=0,1,\ldots, n-1 ,
$$
for some number $\hat C_2$ only dependent on the parameters 
of the class $F_r$ and $P$, and sufficiently large $n$. Hence, 
\be
\int\limits_{t_k}^t \left\| \bar z'_{k}(\xi)- \bar w_{k,s}(\xi) \right\|\,d\xi \leq \hat C_2 \, h_k^{s+2} ,\;\; t\in [t_k, t_{k+1}],
\label{aux2}
\ee
for $s=0,1,\ldots, r-1$.
To bound the second term in (\ref{aux1}), we estimate the difference between two Lagrange polynomials
$$
\left\| \bar w_{k,s} (\xi)-\bar q_{k,s}(\xi)\right\| = 
\Bigl\|
\sum\limits_{p=0}^{s} \prod_{l=0,l\neq p}^{s}\frac{\xi-\xi_{k,l}}{\xi_{k,p}-\xi_{k,l}} \Bigl(\bar f_k(\bar z_{k}(\xi_{k,p}))-\bar f_k(\bar l_{k,s}(\xi_{k,p}))\Bigr)\Bigl\|  
$$
$$
\leq PL\sum_{p=0}^{s}\|\bar z_{k}(\xi_{k,p})-\bar l_{k,s}(\xi_{k,p})\| \prod_{l=0,l\neq p}^{s}\Bigl|\frac{\xi-\xi_{k,l}}{\xi_{k,p}-\xi_{k,l}}\Bigl|
$$
\be
\leq PLC_r\sup\limits_{\xi\in [t_k,t_{k+1}]}\|\bar z_{k}(\xi)-\bar l_{k,s}(\xi)\|,
\label{aux3}
\ee
where $C_r$ is given in (\ref{pom1}). From this we get for $t\in [t_k,t_{k+1}]$
\be
\int\limits_{t_k}^t \left\| \bar w_{k,s} (\xi)-\bar q_{k,s}(\xi)\right\|\, d\xi 
\leq h_k PLC_r\sup\limits_{\xi\in [t_k,t_{k+1}]}\|\bar z_{k}(\xi)-\bar l_{k,s}(\xi)\|.
\label{aux5}
\ee
We now come back to (\ref{aux1}) to get from (\ref{aux2}) and (\ref{aux5}) that 
\be
\sup\limits_{t\in [t_k,t_{k+1}]}\|\bar z_{k}(t)-\bar l_{k,s+1}(t)\| \leq \hat C_2 \, h_k^{s+2} + 
h_k PLC_r\sup\limits_{t\in [t_k,t_{k+1}]}\|\bar z_{k}(t)-\bar l_{k,s}(t)\|,
\label{aux6}
\ee
$s=0,1,\ldots, r-1$, where $\sup\limits_{t\in [t_k,t_{k+1}]} \|\bar z_k(t) - \bar l_{k,0}(t)\leq P(M+LR) h_k$
for sufficiently large $n$. 
By solving the  recurrence inequality with respect to $s$, we obtain  (\ref{Lagr_est_lr}).
\f
We now estimate the global error in $[a,b]$.  We have  
$$
	\| z(t_{k+1})-\bar y_{k+1}\| \leq \| z(t_{k+1})-\bar z_{k}(t_{k+1})\|+\|\bar z_{k}(t_{k+1})-\bar y_{k+1}\| .
$$
Note that $z$ and $\bar z_k$ are solutions of the initial value problems in $[t_k,t_{k+1}]$ with right-hand sides 
$f$ nad $\bar f_k$, and initial conditions $z(t_k)$ and $\bar y_k$, respectively. By a standard use of Gronwall's inequality, 
using (A5) and Lemma 1, we get for $t\in [t_k,t_{k+1}]$ and sufficiently large $n$ 
\be 
\|z(t)-\bar z_k(t)\| \leq e^{Lh_k} \left( \|z(t_k)-\bar y_k\| +h_k \delta(N_k)\right).
\label{aa}
\ee
Hence,
\be
\| z(t_{k+1})-\bar y_{k+1}\| \leq e^{Lh_k} \left( \|z(t_k)-\bar y_k\| +h_k \delta(N_k)\right) +C_1h_k^{r+1},
\label{aa1}
\ee
for $k=0,1,\ldots,n-1$, where  $\| z(t_0)-\bar y_{0}\|=\| \eta-  P_{N_{-1}}\eta\|\leq \gamma(N_{-1})$. 
By solving this recurrence inequality  with respect to $k$, we get
that there is a   number $C_3$  only dependent on the parameters 
of the class $F_r$ and $P$ such that
\begin{equation}
\label{est_bznyn}
	\| z(t_k)-\bar y_{k}\|\leq C_3 \left( \gamma(N_{-1}) +\sum\limits_{j=0}^{k-1} h_j \delta(N_j) +    \sum\limits_{j=0}^{k-1} h_j^{r+1}\right) ,
\end{equation}
for $k=0,1,\ldots,n-1$, for $n$ sufficiently large ($\sum\limits_{j=0}^{-1} =0$). From (\ref{aa}), by slightly changing  the constant $C_3$, we have
for $t\in [t_k, t_{k+1}]$ 
\be
\|z(t)-\bar z_k(t)\| \leq  C_3 \left( \gamma(N_{-1}) +\sum\limits_{j=0}^{k} h_j \delta(N_j) +     \sum\limits_{j=0}^{k-1} h_j^{r+1} \right),
\label{aa2}
\ee
for $k=0,1,\ldots, n-1$, and $n$ sufficiently large. 
\f
By (\ref{Lagr_est_lr}) and (\ref{aa2}) we obtain for $t\in [t_k,t_{k+1}]$ and $k=0,1,\ldots,n-1$ 
$$
	\Bigl\| z(t)-\phi^*_{n,\vN}\Bigl(\mathcal{N}^*_{n,\vN}(f,\eta)\Bigr)(t)\Bigl\| \leq \| z(t)-\bar z_{k}(t)\|+
\Bigl\|\bar z_{k}(t)-\phi^*_{n,\vN}\Bigl(\mathcal{N}^*_{n,\vN}(f,\eta)\Bigr)(t)\Bigl\|
$$
\be
\label{EST_bznfi}
	\leq   C_3 \left( \gamma(N_{-1}) +\sum\limits_{j=0}^{k} h_j \delta(N_j) +     \sum\limits_{j=0}^{k-1} h_j^{r+1} \right)  +C_1h_k^{r+1} ,    
\ee
for $n$ sufficiently large.  This implies   the statement of the theorem in the case $r\geq 1$. 
\f
 For $r=0$, similarly as for $r=1$,  the algorithm  $\phi^*_{n,\vN}$ reduces to the Euler method, i.e., the final approximation is  given by
$$
\bar  l_n(t)= \bar  y_{k} +(t-t_k)\bar f_k(\bar y_{k}), \;  \mbox{ if }\; t\in [t_k,t_{k+1}].
$$
It suffices to note that  the error analysis in  the case $r=1$  only requires  the Lipschitz condition for $\bar f_k$. 
That is, it can be repeated for $r=0$. 
This completes the proof of the theorem.  \qed
\vsn
{\bf Remark 1}$\;\;$  In the special case of uniform discretization of $[a,b]$ and constant truncation parameters, that is,  when $N_{-1}=N_0=…=N_{n-1}=N$, the estimate (\ref{thm1}) can be derived from Lemma 1 above and Theorem 3.3 in \cite{Hein1}. One has to apply the random algorithm from \cite{Hein1}, for a fixed random instant, to the pair $(P_Nf,P_N\eta)$, use Lemma 1 and  note that
information needed for that input is $N$-dimensional. This observation was made by Stefan Heinrich in private communication. 
\vsn
\noindent
 \subsection{ Upper complexity bound }
\noindent
Let $\e>0$. The cost of computing an $\e$-approximation using  information ${\cal N}^*_{n,\vN}$  and the algorithm  $\phi^*_{n,\vN}$
provides an upper bound on the $\e$-complexity. 
The computations    in  $\phi^*_{n,\vN}$  are determined by  vectors  $\vN=[N_{-1},N_0,\ldots, N_{n-1}]$ 
and  $\vM=[M_0,M_1,\ldots, M_{n-1}]$ , where $M_k$ is given by  $M_k=\max\limits_{j=-1,0,\ldots,k} N_j$,  $k=0,1,\ldots, n-1$. 
The number $s$ of information points in ${\cal N}^*_{n,\vN}$ is $O(r^2n)$, with an absolute constant in the $'O'$ notation.
Neglecting the coefficient that only depends on $r$,  the cost of $\phi^*_{n,\vN}$ is thus
\be
\sum\limits_{k=0}^{n-1}  c\left( \max\limits_{j=-1,0,\ldots,k} N_j \right)N_k,
\label{costalg}
\ee
where  $c(N)$ is the cost function defined in Section 2.
Given the mesh points $\{t_k\}$, let ${\rm cost}^*(\e)$ be the 
minimal cost of computing an $\e$-approximation by this class of algorithms,  the minimum
taken with respect to the selection of $n$ and the dimensions $N_{-1},N_0,\ldots, N_{n-1}$.
Let 
\be
U(\e)= \inf \left\{ \sum\limits_{k=0}^{n-1}  c\left( \max\limits_{j=-1,0,\ldots,k} N_j\right) N_k: 
\;\;   \gamma(N_{-1}) +\sum\limits_{j=0}^{n-1} h_j \delta(N_j) +   \sum\limits_{j=0}^{n-1} h_j^{\max\{r,1\} +1}  \leq \e \right\},
\label{mincostalg}
\ee
where the infimum is taken with respect to $n$, $\{h_j\}$ and $\vN$ satisfying the bound. 
For a  given sequence $\{\alpha(n)\}_{n=1}^\infty$, given cost function $c$, and given functions $\gamma $ and $\delta$ defining the class of problems, $U(\e)$ can be computed. 
\f
Due to Theorem 1,   the $\e$-complexity for sufficiently small $\e$  is bounded by 
\be
\comp (\e)\leq {\rm cost}^*(\e)\leq  U(\e/C), 
\label{upp1}
\ee
where $C$ is the constant from Theorem 1.
An obvious choice  is  to take  the truncation parameters  constant in each interval $[t_k,t_{k+1}]$, i.e., $N_{-1}=N_0=\ldots =N_{n-1}=N$.
The minimization in this subclass gives us the value
\be
U^{{\rm eq-dim}} (\e)= \inf \left\{  n  c\left( N\right) N: 
\;\;   \gamma(N) + (b-a) \delta(N) +    \sum\limits_{j=0}^{n-1} h_j^{\max\{r,1\} +1}   \leq \e \right\}. 
\label{mincostalg11}
\ee
We have for sufficiently small $\e>0$ that
\be
\comp (\e) \leq  U(\e/C)  \leq U^{{\rm eq-dim}} (\e/C).
\label{upp2}
\ee
 To further  bound  $U^{{\rm eq-dim}} (\e/C)$ from above,  it suffices to take the value of $ n  c\left( N\right) N$ with the minimal 
$n$ and $N$ such that
$$
C\gamma(N)\leq \e/3 , \;\;\; C(b-a)\delta(N)\leq \e/3,,\;\;\;   C (b-a)\alpha(n)^{\max\{r,1\}}\leq \e/3.
$$
We get 
\vsn
{\bf Proposition 1}$\;\;$ {\it 
There exist positive numbers $C_1$  and $\e_0$   such that for all $\e\in (0,\e_0)$
it holds
\be
 \comp(\e)\leq   n(\e/C_1)\, c(N(\e/C_1) )\,  N(\e/C_1) ,  
\label{prop1}
\ee
where 
$$
n(\e)= \alpha^{-1}\left(\e^{1/\max\{r,1\}} \right) , \;\;
N(\e)= \max\left\{ \gamma^{-1}\left( \e\right) ,   \delta^{-1}\left( \e\right)       \right\}.
$$
}
 \vsn
(For a nonincreasing function $g$ acting from $[1,\infty)$ onto $(0, p]$, ($p>0$), by $g^{-1}$ we mean a function on $(0,p)$
defined by $g^{-1}(\e)=\inf \{ x\in [1,\infty):  g(x)\leq \e \}$.)  
\noindent
{\Large \section{ Lower error and complexity bounds }}
\noindent
In this section we discuss  lower error and complexity bounds. We restrict ourselves to a special (but still interesting) case
 of constant truncation parameters, that is, we assume that $N_{-1}=N_0=\ldots=N_{n-1}=N$. Furthermore, we assume that 
the partitions of $[a,b]$ satisfy the following condition 
(which  most often holds in practice): there exists $\hat K_1$ such that
\be
\alpha(n)\leq \hat K_1n^{-1}, 
\label{al}
\ee
$n=1,2,\ldots \; .$
We shall  denote  information and an algorithm in this case  by ${\cal N}_{n,N}$  and $\phi_{n,N}$, respectively.  
Theorem 1 assures the existence of $C$ (dependent on $\hat K_1$) such that  for $n$ sufficiently large 
\be
e(\phi_{n,N}^* ,{\cal N}_{n,N}^*,  F_r) \leq  C\left(  \gamma(N)+ \delta(N) +n^{-\max\{r,1\}} \right).
\label{th1aa}
\ee
Upper complexity bound of Proposition 1 now holds with 
\be
n(\e)= (1/\e)^{1/\max\{r,1\}},
\label{ne}
\ee
i.e., 
\be
\comp(\e)\leq  (C_1/\e)^{1/\max\{r,1\}}   \, c(N(\e/C_1) )\,  N(\e/C_1) .
\label{prop1aa} 
\ee
\noindent
 \subsection{ Lower error bound }
\noindent
 We now show a matching (up to a constant) lower bound, with respect to  (\ref{th1aa}),
 on  the error of any algorithm $\phi_{n,N}$ using  any information ${\cal N}_{n,N}$
from the considered class. The solution of (\ref{rnie}) for a right-hand side $f$ and
an initial vector $\eta$ will be denoted by $z_{f,\eta}$.  
\vsn
{\bf Theorem 2}$\;\;$ {\it   For any $\hat K$ in (\ref{hatK}) and $\hat K_1$ in (\ref{al}) there exist positive $\hat C$, $n_0$ and $\hat N$ 
such that for any $n\geq n_0$, $N\geq \hat N$, 
for any information ${\cal N}_{n,N}$ and any algorithm $\phi_{n,N}$ it holds
\be
e(\phi_{n,N} ,{\cal N}_{n,N},  F_r) \geq \hat C\left(  \gamma(N)+ \delta(N) +n^{-\max\{r,1\}}\right).
\label{thm2}
\ee
}
{\bf Proof}$\;\;$  Let $(f,\eta), (g,\kappa)\in F_r$ be such that ${\cal N}_{n,N}(f,\eta)={\cal N}_{n,N}(g,\kappa)$. By the triangle inequality we have 
\be
e(\phi_{n,N} , F_r) \geq \frac{1}{2} \sup\limits_{t\in [a,b]}\| z_{f,\eta}(t)- z_{g,\kappa}(t) \|.
\label{p1}
\ee
Using (\ref{rnie}), in a standard way we get that
\be
\sup\limits_{\xi\in [a,t]}\| z_{f,\eta}(\xi)- z_{g,\kappa}(\xi) \| \geq \frac{1}{1+L(t-a)} \left\| \eta-\kappa -\int\limits_{a}^t H(z_{f,\eta}(\xi))\, d\xi \right\|,
\;\; t\in [a,b], 
\label{p2}
\ee
where $H=g-f$.
\f
We now construct suitable pairs $(f,\eta)$ and $(g,\kappa)$. 
\f
Case I.$\;\;$ Let $f=g=0$, $\eta=\gamma(N)e_{N+1}$ and $\kappa=0$. Then $(f,\eta), (g,\kappa)\in F_r$ and ${\cal N}_{n,N}(f,\eta)={\cal N}_{n,N}(g,\kappa)$. Indeed, for example to show (A1), we note that $\|\eta - P_k\eta\| = \gamma(N)\leq \gamma(k)$ for $k\leq N$, and 
$\|\eta - P_k\eta\| = 0$ for $k\geq N+1$.
We have from (\ref{p2})
that 
\be
\sup\limits_{\xi\in [a,b]}\| z_{f,\eta}(\xi)- z_{g,\kappa}(\xi) \| \geq  \gamma(N).
\label{p3}
\ee
Case II.  $\;\;$ Let $f(y)= \delta(N) e_{N+1}$ for $y\in E$ and $N$ sufficiently large so that $\delta(N)\leq M$. Take $g=0$ and $\eta=\kappa$,
where $\eta$ satisfies (A1). 
Then $(f,\eta), (g,\kappa)\in F_r$ and ${\cal N}_{n,N}(f,\eta)={\cal N}_{n,N}(g,\kappa)$. For instance, to see (A5), note that
$\|f(y)-P_kf(y)\|= \delta(N)\leq \delta(k)$ for $k\leq N$, and $\|f(y)-P_kf(y)\|= 0$ for $k\geq N+1$.
From (\ref{p2}) we get
\be
\sup\limits_{\xi\in [a,b]}\| z_{f,\eta}(\xi)- z_{g,\kappa}(\xi) \| \geq  (b-a)\delta(N).
\label{p4}
\ee
Case III. $\;\;$  Let $\eta$ satisfy (A1) and $\kappa=\eta$. We take $f(y)=Ae_1$, where $A>0$. The solution $z_{f,\eta}$ is given by
$$
z_{f,\eta}(t) = A(t-a)e_1 +  \eta,  \;\; t\in [a,b].
$$
Compute the adaptive  information  ${\cal N}_{n,N}(f,\eta)$ for $f$ and $\eta$.
By definition,  ${\cal N}_{n,N}(f,\eta)$ is based on  evaluations of the components, or  partial derivatives  of the components,
 of the function $P_N f$,  at some information points  $\hat y$ such that $\hat y=P_N \hat y$.  The number of these points  is  $O(n)$. 
\f
The function $g$ is defined as $g=f+H$, where $H$ will be given below. Note that the integral in (\ref{p2}) with $t=b$ has now the form
$$
\int\limits_{a}^b H(A(\xi-a)e_1 +  \eta)\, d\xi.  
$$
 Let $r\geq 1$ and $H^{{\rm scal}}:\rr\to \rr$ 
be a nonnegative function such that:
\vsn
$H^{{\rm scal}}\in C^r(\rr)$, $(H^{{\rm scal}})^{(j)}(\hat y^1)=0$ for $j=0,1,\ldots,r$, where $\hat y^1$ is the first component of any information point
(the number of $\hat y^1$ is  $O(n)$),
\f
$H^{{\rm scal}}$ is a Lipschitz function with a constant $L_1$, 
\f
$H^{{\rm scal}}(y^1)\leq M_1$, $\;\;\; | (H^{{\rm scal}})^{(j)}(y^1)|\leq D_1$, for $y^1\in \rr$, $j=1,2, \ldots,r$, for some $M_1$, $D_1$,  and
\be
\int\limits_{a}^b H^{{\rm scal}}(A(\xi-a) +  \eta^1)\, d\xi =\Omega \left(n^{-r}\right).   
\label{integral}
\ee
 For $r=0$ we take the same function $H^{{\rm scal}}$ as for $r=1$.
\f
The construction of such a (bump) function $H^{{\rm scal}}$ is a standard tool when proving lower bounds, see for instance \cite{bk2}. 
We now define for $r\geq 0$ 
$$H(y)= H ^{{\rm scal}} (y^1)\, e_1.$$ 
By taking sufficiently small $A$, $M_1$, $L_1$  and $D_1$, we assure that 
$(f,\eta), (g,\kappa)\in F_r$. Since the derivatives of $H^{{\rm scal}}$ of order $0,1,\ldots,r$ vanish at first component of each information point, 
we have that  
${\cal N}_{n,N}(f,\eta)= {\cal N}_{n,N}(g,\kappa)$. Hence, by     (\ref{p2}) and           (\ref{integral}) we get
\be
\sup\limits_{\xi\in [a,b]}\| z_{f,\eta}(\xi)- z_{g,\kappa}(\xi) \| =\Omega\left(n^{-\max\{r,1\} }\right).
\label{p6}
\ee
The bounds obtained in the three cases above together with (\ref{p1}) lead to the statement of
the theorem. \qed
\vsn
\noindent
 \subsection{ Lower complexity bound  }
\noindent
\vsn
In this section we  discuss a lower $\e$-complexity bound for (\ref{rnie}).
Theorem  2 immediately leads to such a bound  under certain condition, which we     
believe  holds true under  mild assumptions. The condition concerns the number of scalar evaluations $\ell(M_k,N_k)$ in the definition of the complexity.
In our case, we have that $\ell(M_k,N_k)=\ell(N,N)$. The condition  reads:
\vsn
(C) $\;$ for information used by  an  algorithm for solving (\ref{rnie}) with a right-hand side $P_Nf$ and an initial condition $P_N\eta$,
in $\Omega(n)$ time intervals  it holds
 $\ell(N,N) =\Omega(N)$, with  coefficients  in the '$\Omega$' notation only dependent on $\hat K$ and $\hat K_1$ (and the parameters of the class $F_r$, $a$ and $b$).
\vsn
Under  condition (C),  the cost of any algorithm is $\Omega(n\, c(N) N)$.
\vsn 
{\bf Proposition 2 }$\;$ {\it  For any $\hat K$ and $\hat K_1$,  if the class of information  satisfies (C), then there exist positive numbers   
$C_2$ and $\e_0$ such that for all $\e\in (0,\e_0)$ it holds
\be
 \comp(\e)\geq  C_2\, n(\e)\, c(N(\e/C_2) )\, N(\e/C_2) ,  
\label{prop2}
\ee
where 
$n(\e)$ and $N(\e)$ are given in (\ref{ne}) and Proposition 1, respectively.
}
\vsn
{\bf Proof}$\;$ 
Consider an arbitrary algorithm $\phi_{n,N}$ based on some information ${\cal N}_{n,N}$  for which (C) holds. 
 If  $e(\phi_{n,N}, {\cal N}_{n,N}, F_r)\leq \e$, then, due to Theorem 2, we have 
$$
 \hat C\gamma(N)\leq \e , \;\;\; \hat C\delta(N)\leq \e, \;\;\; \hat Cn^{-\max\{r,1\}}\leq \e .
$$
 This yields that 
\be
n\geq \hat C^{1/\max\{r,1\}} n(\e), \;\;\; N\geq N(\e/\hat C).
\label{Nndol}
\ee
Since, by assumption (C), the cost of an algorithm is $\Omega( nc(N)N)$, we get the desired  lower bound. \qed
\vsn
Under condition (C), if $\alpha(n)=O(n^{-1})$,  the lower bound in (\ref{prop2}) matches, up to a constant,  that in  (\ref{prop1aa}). 
\vsn
{\bf Remark 2}$\;\;$ Note that the Taylor algorithm can  potentially be used to solve the finite-dimensional problem in $\rr^N$.
However,  the cost of computing the Taylor information,  which can be as large as $nN^{r+1}c(N)$, 
is much larger than the minimal cost as $N\to \infty$ (unless  function $f$ is very special). 
\noindent
 {\Large\section{ Illustration --  weighted $\ell_p$ spaces}}
\noindent
Consider a countable system of equations
\vsn
\begin{tabular}{ll}
$(z^1)'(t)=f^1(z^1(t),z^2(t),\ldots ),\;\;\;$ & $z^1(a)=\eta^1$\\
$(z^2)'(t)=f^2(z^1(t),z^2(t),\ldots ),\;\;\;$ & $z^2(a)=\eta^2$\\
\vdots                                                    & \vdots\\
\end{tabular}
\vsn
in the interval $[a,b]=[0,1]$.
We embedd this problem into the Banach space setting. Note that
often in applications  components of an infinite sequence $y=(y^1, y^2,\ldots \;)$  are not of the same importance.
Some components may be crucial, while other may even be neglected. It seems reasonable to associate with the components 
 certain positive weights $w_j$, $j=1,2,.\ldots \;\;$.  For $1\leq p<\infty$, we assume that 
$$
 \sum\limits_{j=1}^\infty w_j^p=W^p<\infty.
$$
Consider the Banach space of sequences 
\be
E=\ell_p^w=\left\{ y=(y^1,y^2,\ldots ):\;\; \sum\limits_{j=1}^\infty |y^j|^pw_j^p<\infty\right\},
\label{i2}
\ee
with the norm 
$$
\|y\|_{\ell_p^w}=\left(\sum\limits_{j=1}^\infty  |y^j|^pw_j^p \right)^{1/p}.
$$
The normalized Schauder basis in $\ell_p^w$ is given by
$$
e_j=(0,\ldots, 0, 1/w_j,0,\ldots\; ),
$$
where the $j$th position is nonzero, $j=1,2,\ldots \;$. Note that the basis constant $P$ equals $1$.  
The $j$th coordinate of the sequence $y=(y^1,y^2,\ldots\;)$ in that basis
is given by $w_jy^j$. 
We assume that $\eta=(\eta^1,\eta^2,\ldots)\in \ell_p^w$. The components of the system of equations $f^j$, $j=1,2,\ldots \,$,
 are now treated as functions
$$
f^j:\ell_p^w \to \rr .
$$
We assume about  the components $f^j$ and $\eta^j$ that
\be
|f^j(\eta)|\leq M_j ,\hspace{0.5cm}   |f^j(y)-f^j(\bar y)|\leq L_j\|y-\bar y\|_{\ell_p^w}, \;\; \mbox{ and }\;\; |\eta^j|\leq T_j 
\label{i3}
\ee
for $y,\bar y\in \ell_p^w$, for some nonnegative numbers $M_j$, $L_j$ and $T_j$. 
\f
For $y\in \ell_p^w$, we define $f(y)$ as a sequence 
$$
f(y):=(f^1(y),f^2(y),\ldots ).
$$
We now show, under certain assumptions on $M_j$, $L_j$ and $T_j$, that the conditions (A1)--(A5) hold with $r=0$ for the
pair $(f,\eta)$.
\vsn
{\bf Proposition 3}$\;\;$ {\it 
 Let
\be
W_1:=\left(\sum\limits_{j=1}^\infty M_j^pw_j^p\right)^{1/p}, \;\;  
W_2:=\left(\sum\limits_{j=1}^\infty L_j^pw_j^p\right)^{1/p}
\;\; \mbox{ and } W_3:=\left(\sum\limits_{j=1}^\infty T_j^p w_j^p\right)^{1/p} 
\label{i4}
\ee
be finite numbers. Then
\be
f: \ell_p^w \to \ell_p^w,
\label{i5}
\ee
\be
\|f(\eta)\|_{\ell_p^w} \leq W_1,
\label{i6}
\ee
\be
\|f(y)-f(\bar y)\|_{\ell_p^w} \leq W_2 \|y-\bar y\|_{\ell_p^w},
\label{i7}
\ee
\be 
\|f(y)-P_kf(y)\|_{\ell_p^w} \leq 2^{1-1/p} \left(\sum\limits_{j=k+1}^\infty M_j^pw_j^p +\|y-\eta\|_{\ell_p^w}^p \sum\limits_{j=k+1}^\infty L_j^pw_j^p\right)^{1/p},
\label{i8}
\ee
and
\be
\|\eta -P_k \eta\|_{\ell_p^w} \leq \left( \sum\limits_{j=k+1}^\infty T_j^p w_j^p \right)^{1/p}.
\label{i9}
\ee
}
{\bf Proof}$\;\;$ Note that
\be
|f^j(y)|^p\leq \left( |f^j(\eta)|+L_j\|y-\eta\|_{\ell_p^w} \right)^p\leq 2^{p-1}\left( |f^j(\eta)|^p + L_j^p\|y-\eta\|_{\ell_p^w}^p\right),
\label{i10}
\ee
which gives   
$$
\| f(y)\|_{\ell_p^w}^p\leq  2^{p-1} \left( W_1^p + W_2^p \|y-\eta\|_{\ell_p^w}^p\right).
$$
Hence, (\ref{i5}) holds. The proofs of (\ref{i6}) and (\ref{i7}) are straightforward. 
\f
To show (\ref{i8}), we see that
$$
\|f(y)-P_kf(y)\|_{\ell_p^w}^p  = \sum\limits_{j=k+1}^\infty |f^j (y)|^p w_j^p ,
$$
and we use (\ref{i10}). We get 
$$
\|f(y)-P_kf(y)\|_{\ell_p^w}^p \leq  \sum\limits_{j=k+1}^\infty 2^{p-1}\left( |f^j(\eta)|^p w_j^p + L_j^p\|y-\eta\|_{\ell_p^w}^p w_j^p \right),
$$
which yields (\ref{i8}). The proof of (\ref{i9}) is similar.  
\qed
\vsn
Proposition 3 yields that the system that we started with has the form 
$$
z'(t)=f(z(t)), \;\;\; t\in [0,1], \;\;\; z(0)=\eta
$$
in the Banach space $E=\ell_p^w$, and the pair  $(f,\eta)$ satisfies (A1)--(A5) with $r=0$. The sequences $\gamma(k)$ and $\delta(k)$ are defined 
by the right-hand sides of (\ref{i9}) and (\ref{i8}), respectively, taking into account that $y\in K(\eta,R)$. 
\f
There exists a unique solution $z:[a,b]\to \ell_p^w$,
which can be approximated using the truncated Euler algorithm described in Section 3. 
\f
We apply  the Euler algorithm  on the uniform mesh with $h_k=(b-a)/n$ and $N_{-1}=N_0=\ldots=N_{n-1}=N$.
The Euler approximation in $[a,b]$ to the solution $z$, where 
$$
z(t)=(z^1(t), z^2(t),\ldots,  z^N(t),  z^{N+1}(t), \ldots \;\;) \in \ell_p^w
$$ 
is defined by  
$$
\bar l_n(t)=( \bar l_n^1(t), \bar l_n^2(t),\ldots,  \bar l_n^N(t), 0,0,\ldots \;)
$$
Consider the following  example. Let $p\in (1,\infty)$, $w_j=1/j$, $M_j=L_j=T_j=1$. Then $W_1=W_2=W_3=W=
\left(\sum\limits_{j=1}^\infty (1/j)^p\right)^{1/p}$.
Furthermore, $M=L=W$ and 
$R$ of Lemma 1 are  known numbers.  
Since 
$$
\sum\limits_{j=k+1}^\infty (1/j)^p \leq \int\limits_k^\infty (1/x)^p\, dx=(1/k)^{p-1}/(p-1),
$$
we can take $\gamma(k)=\delta(k)/(2R)= (1/k)^{1-1/p}/(p-1)^{1/p}$. We also have $W\leq (p/(p-1))^{1/p}$.
One can see that  the error of the truncated Euler algorithm in the class $F_0$ is bounded by 
$$
e(\phi^*_{n,N}, {\cal N}^*_{n,N}, F_0)  \leq A   (1/N)^{1-1/p} + B(1/n),
$$
where 
$$
A=  (2R+1)\exp(W)/(p-1)^{1/p}, \;\;\; B=  W(R+1)(3\exp(W) -2)
$$
 are known absolute constants.
The cost of the truncated Euler algorithm is equal to $nc(N)N$. 
\f
Let $\e>0$.  Take, for instance, the cost function $c(N)=N^{\beta}$, $\beta\geq 0$, and consider 
 minimization of the cost of the algorithm with the error bounded by $\e$:
\vsn
minimize $nc(N)N=nN^{1+\beta}$ $\;\;\;\;$ subject to  $\;\;\;A   (1/N)^{1-1/p} + B(1/n)\leq \e$.
\vsn
That is,   we wish to find the best (in the framework of the example) discretization and truncation parameters $n$ and $N$.
The solution is given by 
$$
 n=n(\varepsilon)=\Bigl\lceil\frac{B(p(\beta+2)-1)}{(p-1)\varepsilon}\Bigr\rceil
$$
and 
$$
N=N(\varepsilon)=\Bigl\lceil\frac{A(\beta+2-1/p)}{(\beta+1)\varepsilon}\Bigr\rceil^{p/(p-1)}.
$$ 
The minimal cost is then equal to 
$$
n(\varepsilon)(N(\varepsilon))^{1+\beta}=O\Biggl(\Bigl(\frac{1}{\varepsilon}\Bigr)^{(p(\beta+2)-1)/(p-1)}\Biggr).
$$
For example, if $p=2$ and $\beta=1$, then the minimal cost is $O((1/\varepsilon)^{5})$. The constants in the $'O'$ notation 
are known absolute constants. For comparison, if we solve a finite-dimensional system of $N$ equations with fixed $N$, then the cost
of the Euler algorithm is  $O(1/\varepsilon)$, and $N$ enters the constant.
\vsn
{\bf Remark 3} $\;$
Note that the weights $w_j$ may enter the formulation of the problem through  the norm in $E$, or through the choice of a subclass
of problems to be considered. Let, for example,  $E=\ell_p^w$ with the norm defined above. Consider the class of 
right-hand side functions $f:\ell_p^w \to \ell_p^w$ and initial conditions defined by (\ref{i3}) with $M_j=L_j=T_j=1$, and the problem
$$
z'(t)=f(z(t)), \;\; z(0)=\eta,\;\;\; t\in [0,1].
$$
Here the weights define the norm of the space, and they do not enter the definition of the class of right-hand side functions
and initial conditions.
\f
The problem can reformulated as follows. For a sequence $y=(y^1, y^2, \ldots\;\;)\in \ell_p$, let 
$\tilde f(y)$ be a sequence $(\tilde f^1(y),\tilde f^2(y),\ldots\;\;)$, where
$$
\tilde f^j(y)= w_j f^j(y^1/w_1, y^2/w_2, \ldots\;\;).
$$
Then $\tilde f:\ell_p\to \ell_p$ with the standard norm (which does not depend on $w_j$). Consider the problem
$$
\tilde z '(t)=\tilde f (\tilde z (t)), \;\; \tilde z (0)=\tilde \eta,\;\;\; t\in [0,1],
$$
with $\tilde \eta =(w_1\eta^1, w_2\eta^2,\ldots \;\;)$. Note that both initial value problems are equivalent, since
$$
\tilde z^j(t)= w_j z^j(t) \;\; \mbox{ for } j=1,2\ldots\;\;.
$$
The restrictions (\ref{i3}) are equivalent to the following restrictions on $\tilde f^j$ and $\tilde\eta^j$
$$
\frac{1}{w_j} |\tilde f^j(\tilde \eta)| \leq 1, \;\;\;\; \frac{1}{w_j} |\tilde f^j(y)-\tilde f^j(\bar y)|\leq \|y-\bar y\|_{\ell_p}, \;\;\;\; \frac{1}{w_j}| \tilde\eta^j|\leq 1.
$$
In the alternative  formulation, the weights  appear in the restrictions on $\tilde f$ and $\tilde\eta$, not in the norm of the space, 
see e.g. \cite{deiml}, p. 109.
\vsn
\noindent
{\Large \section{ Auxilliary facts}}  
\noindent
For  convenience of the reader, we recall some well known facts that are used in this paper.
Let $\alpha:[a,b]\to E$. We define the Riemann integral 
$$
\int\limits_a^b \a (t)\, dt   
$$
(an element of $E$) in the same way as we do for real functions as a limit of Riemann sums, see e.g. \cite{LusSob} or  \cite{Rall}. If 
$\a (t)=\sum\limits_{j=1}^\infty \a ^j(t)e_j$ is Riemann integrable, then $\a^j$ are also  Riemann integrable real functions and  
$$
\int\limits_a^b \a (t)\, dt = \sum\limits_{j=1}^\infty \left( \int\limits_a^b \a ^j (t)\, dt \right)\, e_j.
$$
We have that $\a_1(t)\, e_1= P_1\a(t)$ and $\a ^j(t)\, e_j= (P_j -P_{j-1}) \a(t)$ for $j\geq 2$.
If $\a$ is a continuous function then $\a ^j$ are continuous, if $\a$ is a Lipschitz function with a constant $C$, then
$\a ^j$ are Lipschitz functions with the constant $2PC$. 
\f
Let $\a$ be $k$ times Frech\'et differentiable in $[a,b]$. The derivative  $\a^{(k)}(t)$ can be 
identified with an element of $E$, i.e., $\a^{(k)}: [a,b]\to E$. Then $\a^j$ are also $k$ times differentiable functions, and 
$$
\a ^{(k)} (t)=\sum\limits_{j=1}^\infty (\a ^j)^{(k)} (t)e_j.
$$
Let $a=t_0<\ldots<t_p=b$. Let $w:[a,b]\to E$ be an interpolation  polynomial of degree at most $p$ 
(i.e., a function of the form $w(t)=\sum\limits_{i=0}^p t^i a_i $ for some $a_i\in E$) such that
$$
w(t_j)=\a(t_j), \;\;\; j=0,1,\ldots, p.
$$
As in the case of real-valued functions, one can see that the interpolation conditions are satisfied for 
$$
w(t)=\sum\limits_{i=0}^p \prod\limits_{s=0,s\ne i}^p \frac{t-t_s}{t_i-t_s}\, \a(t_i).
$$
The  coefficients $w^j$ of $w$ in the basis $\{e_j\}$ are real-valued polynomials of degree at most $p$. 
From  the interpolation conditions for $w^j$, we see that  $w^j$ are unique, and consequently  so is $w$. 
 \f
Let $\a\in C^{p+1}([a,b],E)$.    In the same way as for  real-valued 
functions, one can prove that the remainder $R(t)=\a(t)-w(t)$ of the interpolation formula can be written in the integral form as
$$
R(t)=\prod\limits_{i=0}^p (t-t_i)\, G(t),
$$ 
where $G(t)\in E$ is given by (see e.g. \cite{jank})
\be
G(t)= \int\limits_0^{1}\int\limits_0^{\xi_0}\ldots\int\limits_0^{\xi_{p-1}} \a^{(p+1)}\left( t+\xi_0(t_0-t)+\ldots +\xi_p(t_p-t_{p-1} ) \right)\,
d\xi_p\ldots d\xi_1d\xi_0.
\label{interp-remain}
\ee
Let $D\subset E$ be an open nonempty subset, and  $g:D\to E$ be a $k$ times Frech\'et differentiable function in $D$,
$$
g(y)=\sum\limits_{j=1}^\infty g^j(y) e_j.
$$
Since $g^1(y)e_1=P_1g(y)$ and $g^j(y)e_j= (P_j-P_{j-1})g(y)$ for $j\geq 2$,  by the definition of Frech\'et derivative we see that
$g^j$ are also $k$ times Frech\'et differentiable functions. If $\|g^{(k)}(y)\| \leq Z$ for some constant $Z$, then
$\| (g^j)^{(k)}(y) \| \leq 2PZ$, where the first symbol $\|\cdot\|$ means the norm of a $k$-linear operator in $E^k$,  while 
the second one means the norm of a $k$-linear functional.  
\vsn
\section{Conclusions}
\vsn
We analyzed the finite-dimensional solution of inital value problems in infinite-dimensional Banach spaces. For $r$-smooth right-hand side functions, we showed an algorithm 
for solving such problems on a non-uniform mesh with variable dimensions. For a constant  dimension $N$, under additional assumptions,
we proved  its error and cost optimality (up to constants), as the truncation and discretization parameters $N$ and $n$ tend to infinity.
The results were illustrated by a countable system in the weighted $\ell_p$ space.
\vsn
{\bf Acknowledgments}$\;\;$ We are indebted to Stefan Heinrich for reading the manuscript and giving his comments, see Remark 1.
\vsn

\end{document}